\documentclass{amsart}
\usepackage{amsmath,amssymb,amsthm}
\usepackage[all]{xy}

\newcommand{\bs}{\backslash}
\newcommand{\longto}{\longrightarrow}

\newcommand{\R}{\mathbb{R}}

\newcommand{\Hom}{\mathrm{Hom}}
\newcommand{\Rep}{\mathrm{Rep}}

\newcommand{\calC}{\mathcal{C}}
\newcommand{\piS}{\pi_1(S^2\bs\{s_1,\,\ldots\, ,s_l\})}

\newcommand{\w}{\omega}

\newcommand{\im}{\mathrm{Im~}}
\newcommand{\pconj}{\calC_1\times\cdots\times\calC_l}
\newcommand{\calD}{\mathcal{D}}

\newcommand{\wred}{\w^{red}}
\newcommand{\quot}{//}

\newcommand{\calN}{\mathcal{N}}
\newcommand{\calB}{\mathcal{B}}

\theoremstyle{plain}

\newtheorem{thm}{Theorem}[section]
\newtheorem{prop}[thm]{Proposition}
\newtheorem{cor}[thm]{Corollary}
\newtheorem{lem}[thm]{Lemma}

\theoremstyle{definition}
\newtheorem{defi}[thm]{Definition}

\newtheorem{rk}[thm]{Remark}
\newtheorem*{ack}{Acknowledgements}

\title[Quasi-Hamiltonian quotients as disjoint unions of symplectic manifolds]{Quasi-Hamiltonian quotients as disjoint unions\\ of symplectic manifolds}
\author{Florent Schaffhauser - Keio University, Yokohama, Japan}
\date{}

\begin{document}

\begin{abstract}
The main result of this paper is Theorem \ref{reduction_strat} which says that the quotient $\mu^{-1}(\{1\})/U$ associated to a quasi-Hamiltonian space $(M,\w,\mu:M\to U)$ has a symplectic structure \emph{even when} $1$ \emph{is not a regular value of the momentum map} $\mu$. Namely, it is a disjoint union of symplectic manifolds of possibly different dimensions, which generalizes the result of Alekseev, Malkin and Meinrenken in \cite{AMM}. We illustrate this theorem with the example of representation spaces of surface groups. As an intermediary step, we give a new class of examples of quasi-Hamiltonian spaces: the isotropy submanifold $M_K$ whose points are the points of $M$ with isotropy group $K\subset U$.
\end{abstract}

\keywords{Momentum maps, quasi-Hamiltonian spaces, moduli spaces}

\maketitle

The notion of quasi-Hamiltonian space was introduced by Alekseev, Malkin and Meinrenken in their paper \cite{AMM}. The main motivation for it was the existence, under some regularity assumptions, of a symplectic structure on the associated quasi-Hamiltonian quotient. Throughout their paper, the analogy with usual Hamiltonian spaces is often used as a guiding principle, replacing Lie-algebra-valued momentum maps with Lie-group-valued momentum maps. In the Hamiltonian setting, when the usual regularity assumptions on the group action or the momentum map are dropped, Lerman and Sjamaar showed in \cite{LS} that the quotient associated to a Hamiltonian space carries a stratified symplectic structure. In particular, this quotient space is a disjoint union of symplectic manifolds.. In this paper, we prove an analogous result for quasi-Hamiltonian quotients. More precisely, we show that for any quasi-Hamiltonian space $(M,\w,\mu:M\to U)$, the associated quotient $M//U:=\mu^{-1}(\{1\})/U$ is a disjoint union of symplectic manifolds (Theorem \ref{reduction_strat}): $$\mu^{-1}(\{1\})/U=\bigsqcup_{j\in J} (\mu^{-1}(\{1\})\cap M_{K_j})/ L_{K_j}.$$ Here $K_j$ denotes a closed subgroup of $U$ and $M_{K_j}$ denotes the isotropy submanifold of type $K_j$: $M_{K_j}=\{x\in M\ |\ U_x=K_j\}$. Finally, $L_{K_j}$ is the quotient group $L_{K_j}=\mathcal{N}(K_j)/K_j$, where $\mathcal{N}(K_j)$ is the normalizer of $K_j$ in $U$. As an intermediary step in our study, we show that $M_{K_j}$ is a quasi-Hamiltonian space when endowed with the (free) action of $L_{K_j}$.

\begin{ack}
This paper was written during my post-doctoral stay at Keio University, which was made possible thanks to the support of the Japanese Society for Promotion of Science (JSPS). I would like to thank the referee for comments and suggestions to improve the paper.
\end{ack}

\section{Quasi-Hamiltonian spaces}

\subsection{Definition}

Throughout this paper, we shall designate by $U$ a compact connected Lie group whose Lie algebra
$\mathfrak{u}=Lie(U)=T_1 U$ is equipped with an $Ad$-invariant
positive definite product denoted by $(.\,|\,.)$. We denote by
$\chi$ (half) the Cartan $3$-form of $U$, that is, the left-invariant
$3$-form on $U$ defined on $\mathfrak{u}=T_1 U$ by: $$\chi_1(X,Y,Z) :=
\frac{1}{2} (X \,|\, [Y,Z]) = \frac{1}{2} ([X,Y] \,|\, Z).$$
Recall that, since $(.\,|\, .)$
is $Ad$-invariant, $\chi$ is also right-invariant and that it is a closed form. Further, let us denote by
$\theta^L$ and $\theta^R$ the respectively left-invariant and right-invariant Maurer-Cartan
$1$-forms on $U$: they take value in $\mathfrak{u}$ and are the identity on
$\mathfrak{u}$, meaning that for any $u \in U$ and any $\xi\in T_u U$,
$$\theta^L_u(\xi)=u^{-1}.\xi \mathrm{\quad and\quad } \theta^R_u(\xi)=\xi.u^{-1}$$
 (where we denote by a point $.$ the effect of translations on
tangent vectors). Finally, we denote by $M$ a manifold on which the group
$U$ acts, and by $X^{\#}$ the fundamental vector field on $M$ defined,
for any $X\in\mathfrak{u}$, by the action of $U$ in the following way:
\begin{eqnarray*}
& X^{\#}_x := \frac{d}{dt}|_{t=0}\
\big(\exp(tX).x\big) & 
\end{eqnarray*}
 for any $x\in M$. We then recall the definition of a quasi-Hamiltonian space, which was first introduced in \cite{AMM}.
\begin{defi}[Quasi-Hamiltonian space, \cite{AMM}]\label{def-qHam}
Let $(M,\w)$ be a manifold endowed with a $2$-form $\w$ and an action of the
Lie group $(U,(.\,|\, .))$ leaving the $2$-form $\w$ invariant. Let $\mu:M
\to U$ be a $U$-equivariant map (for the conjugacy action of $U$ on itself).\\ Then
$(M,\w,\mu:M\to U)$ is said to be a \emph{quasi-Hamiltonian space} with respect to the
action of $U$ if the map $\mu: M\to U$ satisfies the following three
conditions:
\begin{enumerate}
\item[(i)] $d\w = -\mu^* \chi$
\item[(ii)] for all $x \in M$, $\ker \w_x= \{X^{\#}_x\ :\
X\in\mathfrak{u}\,|\, (Ad\,\mu(x)+Id).X=0\}$
\item[(iii)] for all $X\in \mathfrak{u}$, $\iota_{X^{\#}}\w =
\frac{1}{2}\mu^*(\theta^L +\theta^R \,|\, X)$
\end{enumerate}
 where $(\theta^L+\theta^R\,|\, X)$ is the real-valued $1$-form
defined on $U$ for any $X\in\mathfrak{u}$ by $(\theta^L + \theta ^R \,|\,
X)_u(\xi) := (\theta^L_u(\xi) + \theta^R_u(\xi)\,|\, X)$ (where $u\in U$ and
$\xi \in T_u U$).\\ In analogy with the usual Hamiltonian case, the map
$\mu$ is called the \emph{momentum map}. 
\end{defi}

\subsection{Examples}

In this subsection, we recall the fundamental examples of quasi-Hamiltonian spaces. We will use them in section \ref{sympstructure} to illustrate Theorem \ref{reduction_strat}.
\begin{prop}[\cite{AMM}]\label{conjclass}
Let $\calC \subset U$ be a conjugacy class of a Lie group $(U,(.\,|\, .))$.
The tangent space to $\calC$ at $u\in \calC$ is $T_u\calC =\{ X.u-u.X\ : \
X\in\mathfrak{u}\}$. The $2$-form $\w$ on $\calC$ given at $u\in\calC$ by
$$\w_u(X.u - u.X,Y.u -u.Y) = \frac{1}{2}\big((Ad\, u.X \,|\, Y) - (Ad\,
u.Y\,|\, X)\big)$$  is well-defined and makes $\calC$ a
quasi-Hamiltonian space for the conjugacy action with momentum map the
inclusion $\mu:\calC \hookrightarrow U$. Such a $2$-form is actually unique.
\end{prop}
The following theorem explains how to construct a new quasi-Hamiltonian $U$-space out of two existing quasi-Hamiltonian $U$-spaces.
\begin{thm}[Fusion product of quasi-Hamiltonian spaces, \cite{AMM}]\label{fusion}
Let $(M_1,\w_1,\mu_1)$ and $(M_2,$ $\w_2,\mu_2)$ be two quasi-Hamiltonian
$U$-spaces. Endow $M_1 \times M_2$ with the diagonal action of $U$. Then the
$2$-form $$\w:=(\w_1\oplus \w_2) + \frac{1}{2}(\mu_1^*\theta^L \wedge
\mu_2^*\theta^R)$$ makes $M_1\times M_2$ a quasi-Hamiltonian space with
momentum map:
$$\begin{array}{rrcl}
\mu_1\cdot\mu_2: & M_1 \times M_2 & \longto & U \\
& (x_1,x_2) & \longmapsto & \mu_1(x_1)\mu_2(x_2)
\end{array}$$
\end{thm}
\begin{cor}\label{pconj}
The product $\pconj$ of $l$
conjugacy classes of $U$ is a quasi-Hamiltonian space for the diagonal
action of $U$, with momentum map the product $\mu(u_1,\,\ldots\, ,u_l)=
u_1\ldots u_l$. 
\end{cor}
\begin{prop}[\cite{AMM}]\label{intfuseddouble}
The manifold $\mathfrak{D}(U):=U\times U$ equipped with the diagonal
conjugacy action of $U$, the $U$-invariant $2$-form $$\w =
\frac{1}{2}(\alpha^*\theta^L\wedge\beta^*\theta^R) +
\frac{1}{2}(\alpha^*\theta^R\wedge\beta^*\theta^L) +
\frac{1}{2} \big((\alpha\cdot\beta)^*\theta^L\wedge
(\alpha^{-1}\cdot\beta^{-1})^*\theta^R\big)$$ and the
equivariant momentum map
$$\begin{array}{rrcl}
\mu : & \mathfrak{D}(U)=U\times U & \longto & U \\
& (a,b) & \longmapsto & aba^{-1}b^{-1}
\end{array}$$
 (where $\alpha$ and $\beta$ are the projections respectively on the first and
second factors of $\mathfrak{D}(U)$) is a quasi-Hamiltonian $U$-space, called
the \emph{internally fused double of} $U$.
\end{prop}
\begin{cor}
The product manifold 
$$\mathcal{M}_{g,l}:=\underset{g\
\mathrm{times}}{\underbrace{(U\times U)\times\cdots\times(U\times U)}}\times\pconj$$
equipped with the diagonal $U$-action and the momentum map
$$\begin{array}{rrcl}
\mu_{g,l}: & (U\times U)\times\cdots\times(U\times U)\times\pconj &
\longto & U \\
& (a_1,b_1,\, {\ldots} \, ,a_g,b_g,u_1, \, {\ldots} \, ,u_l) & \longmapsto &
[a_1,b_1]{\ldots}[a_g,b_g]u_1{\ldots} u_l 
\end{array}$$
is a quasi-Hamiltonian space.
\end{cor}
This space plays a very important role in the description of
symplectic structures on representation spaces of fundamental groups of
Riemann surfaces (see \cite{AMM} and section \ref{sympstructure} below).

\subsection{Properties of quasi-Hamiltonian spaces}

We now give the properties of quasi-Hamiltonian spaces that we shall need
when considering the reduction theory of quasi-Hamiltonian spaces. The results in the Proposition below are quasi-Hamiltonian analogues of classical lemmas entering the reduction theory for usual Hamiltonian
spaces.
\begin{prop}[\cite{AMM}]\label{facts}
Let $(M,\w,\mu:M\to U)$ be a quasi-Hamiltonian $U$-space and let $x\in M$.
Then:
\begin{enumerate}
\item[(i)] The map $$\begin{array}{rrcl}
\Lambda_x : & \ker (Ad\, \mu(x) + Id) & \longto & \ker \w_x \\
& X & \longmapsto & X^{\#}_x = \frac{d}{dt}|_{t=0}\ \big(\exp(tX).x\big)
\end{array}$$  is an isomorphism.
\item[(ii)] $\ker T_x\mu \cap \ker \w_x = \{0\}$
\item[(iii)] The left translation $$\begin{array}{rcl}
U & \longto & U \\
u & \longmapsto & \big(\mu(x)\big)^{-1}u
\end{array}$$  induces an isomorphism $$\im T_x\mu \simeq
\mathfrak{u}_x^{\perp}$$ where $\mathfrak{u}_x=\{X\in \mathfrak{u} \ |\
X^{\#}_x=0\}$ is the Lie algebra of the stabilizer $U_x$ of $x$ and
$\mathfrak{u}_x^{\perp}$ denotes its orthogonal with respect to $(.\,|\,
.)$. Equivalently, $\im (\mu^*\theta^L)_x=\mathfrak{u}_x^{\perp}$ (and
likewise, $\im (\mu^*\theta^R)_x = \mathfrak{u}_x^{\perp}$). 
\item[(iv)] $(\ker T_x\mu)^{\perp_{\w}} = \{X^{\#}_x \, :\,X\in
\mathfrak{u}\}$, where $(\ker T_x\mu)^{\perp_{\w}} \subset T_x M$ denotes
the subspace of $T_x M$ orthogonal to $\ker T_x\mu$ with respect to $\w_x$.
\end{enumerate}
\end{prop}
We end this subsection with a result that we will need in subsection \ref{reduction}. This theorem relates quasi-Hamiltonian spaces to usual Hamiltonian spaces and we quote it from \cite{AMM} (see remark 3.3, see also \cite{HJS}).
\begin{thm}[Linearization of quasi-Hamiltonian spaces, \cite{AMM}]\label{linearization}
Let $(M_0,\w_0,\mu_0:M_0\to U)$ be a quasi-Hamiltonian $U$-space. Suppose there exists an $Ad$-stable open subset $\calD\subset \mathfrak{u}$ such that $\exp|_{\calD}:\calD\to \exp(\calD)$ is a diffeomorphism onto a open subset of $U$ \emph{containing} $\mu_0(M_0)$. Denote by $\exp^{-1}:\exp(\calD)\to \calD$ the inverse of $\exp|_{\calD}$. Then, there exists a \emph{symplectic} $2$-form $\widetilde{\w_0}$ on $M_0$ such that $(M_0,\widetilde{\w_0},\widetilde{\mu_0}:=\exp^{-1}\circ\mu_0:M_0\to\mathfrak{u})$ is a Hamiltonian $U$-space in the usual sense, for the \emph{same} $U$-action. Furthermore, one has: $$\mu_0^{-1}(\{1_U\})=\widetilde{\mu_0}^{-1}(\{0\})$$ and therefore
$$\mu_0^{-1}(\{1_U\})/U=\widetilde{\mu_0}^{-1}(\{0\})/U$$
\end{thm}

\section{Reduction theory of quasi-Hamiltonian spaces}\label{reductiontheory}

In this section. we show that for any quasi-Hamiltonian space $(M,\w,\mu:M\to U)$, the associated quotient $M^{red}:=\mu^{-1}(\{1\})/U$ is a disjoint union of symplectic manifolds. We begin by reviewing the usual Hamiltonian case and the reduction theorem of Alekseev, Malkin and Meinrenken for quasi-Hamiltonian spaces (Theorem \ref{reduction}). We then begin our study of the stratified case and prove the main result of this paper (Theorem \ref{reduction_strat}). We also prove that isotropy submanifolds are always quasi-Hamiltonian spaces (Theorem \ref{strata_one}).

\subsection{Symplectic reduction in the usual Hamiltonian setting}

In this subsection, we recall how to obtain a
symplectic manifold from a usual Hamiltonian space by a reduction procedure,
that is to say, by taking the quotient of a fiber $\mu^{-1}(\{u\})$ of the
momentum map by the action of the stabilizer group $U_u$, which preserves
the fiber $\mu^{-1}(\{u\})$ since $\mu$ is equivariant. This reduction procedure is usually called the Marsden-Meyer-Weinstein procedure.\vskip 9pt

Let us first recall how to obtain differential forms on an orbit space $N/G$
where $N$ is a manifold acted on by a Lie group $G$. We will assume that $G$ is compact and that it acts freely on $N$ so that $N/G$ is a manifold and the submersion $p: N\to N/G$ is
a locally trivial principal fibration with structural group $G$. Let $[x]$
denote the $G$-orbit of $x \in N$. Since $p$ is surjective, one has $T_{[x]}(N/G)=\mathrm{Im~} T_x p\simeq T_xN/\ker T_xp$. And $\ker T_xp$
consists exactly of the vectors tangent to $N$ at $x$ which are actually
tangent to the $G$-orbit of $x$ in $N$. Those are exactly the values at $x$
of fundamental vector fields: $$\ker T_x p= T_x(G.x)=\{X^{\#}_x\ :\
X\in\mathfrak{g}=Lie(G)\}.$$ Let then $\alpha$ be a differential form on $N$
(say, a $2$-form). Under what conditions does $\alpha$ define a $2$-form
$\overline{\alpha}$ on $N/G$ verifying $p^*\overline{\alpha}=\alpha$\, ?
This last condition amounts to saying that
$\overline{\alpha}_{[x]}([v],[w])=\alpha_x(v,w)$ for all $x\in N$ and all
$v,w \in T_x N$. One then checks that the left-hand side term of this
equation is well-defined by this relation if and only if the $2$-form
$\alpha$ is $G$-invariant. Further, since $X^{\#}_x$ is sent to $0$ in
$T_{[x]}(N/G)$ by the map $T_x p$, the relation $p^*\overline{\alpha}=\alpha$ implies that
$\iota_{X^{\#}}\alpha = 0$ for all $X\in\mathfrak{g}$. These two conditions turn
out to be enough:
\begin{lem}\label{basicforms}
Let $p: N\to B=N/G$ be a locally trivial principal fibration with
structural group $G$ and let $\alpha$ be a differential form on $N$. If
$\alpha$ satisfies
\begin{eqnarray*}
& g^*\alpha = \alpha & \mathrm{for\ all\ } g\in G\quad
\mathrm{(}G\mathrm{-invariance)}\\
\mathrm{and} & & \\
& \iota_{X^{\#}}\alpha=0 & \mathrm{for\ all\ } X\in \mathfrak{g}=Lie(G) 
\end{eqnarray*}
 then there exists a unique differential form $\overline{\alpha}$ on $B$
satisfying $p^*\overline{\alpha}=\alpha$. In such a case, the differential form $\alpha$
on $N$ is said to be \emph{basic}.
\end{lem}
 Observe that if $G$ is compact and connected (so that the
exponential map is surjective), the condition $g^*\alpha=\alpha$ for all
$g\in G$ may be replaced by $\mathcal{L}_{X^{\#}}\alpha=0$ for all $X\in
\mathfrak{g}$ (which is always implied by the $G$-invariance). Further, observe that if $\alpha$ is basic then $d\alpha$ is also basic (the first condition is obvious and the second follows from the Cartan homotopy formula: $\iota_{X^{\#}}(d\alpha)=\mathcal{L}_{X^{\#}}\alpha - d(\iota_{X^{\#}}\alpha)$).\vskip 9pt

We can now use this result to construct differential forms on orbit spaces
associated to level manifolds of the momentum map. Let $(M,\w)$ be a symplectic
manifold endowed with a Hamiltonian action of a compact connected Lie group $U$ with momentum
map $\mu:M\to\mathfrak{u}^*$, and take $N:=\mu^{-1}(\{\zeta\})$ where $\zeta\in
\mathfrak{u}^*$. Because of the equivariance of $\mu$, the stabilizer
$G:=U_{\zeta}$ of $\zeta$ for the co-adjoint action of $U$ on $\mathfrak{u}^*$
acts on $N=\mu^{-1}(\{\zeta\})$. Assuming that $U_{\zeta}$ (which is compact) acts freely on $\mu^{-1}(\{\zeta\})$, one has that $\zeta$ is a regular value of
$\mu$ (see the proof of
Theorem \ref{reduction} for similar reasoning) and we then have a
principal fibre bundle $p:\mu^{-1}(\{\zeta\})\to\mu^{-1}(\{\zeta\})/U_{\zeta}$ and
the following diagram:
\begin{center}$\xymatrix{
{\mu^{-1}(\{\zeta\})}\ \ar[d]_p \ar@{^{(}->}[r]^i & M\\
{\mu^{-1}(\{\zeta\})/U_{\zeta}}
}$\end{center}
 where $i:\mu^{-1}(\{\zeta\})\hookrightarrow M$ is the
inclusion map. The $2$-form $\w$ on $M$ induces a $2$-form $i^*\w$ on
$\mu^{-1}(\{\zeta\})$, which turns out to be basic (again, see the proof of
Theorem \ref{reduction} for similar reasoning). Therefore, by Lemma
\ref{basicforms}, there exists a unique
$2$-form $\wred$ on $\mu^{-1}(\{\zeta\})/U_{\zeta}$ such that $p^*\wred =
i^*\w$. Since $\w$ is closed, so is $\wred$. And one may then notice that a
vector $v\in T_x N=\ker T_x\mu$ is sent by $T_x p$ to a vector in $\ker \wred_{[x]}$ if
and only if $v$ is contained in $(T_x N)^{\perp_{\w}} = (\ker T_x\mu)^{\perp_{\w}} =
\{X^{\#}_x\, :\, X\in \mathfrak{u}\}$ as well. But then $v=X^{\#}_x \in \ker
T_x\mu\cap (\ker T_x\mu)^{\perp_{\w}}$, so that by the equivariance of
$\mu$, one has, denoting by $X^{\dagger}$ the fundamental vector field on
$\mathfrak{u}^*$ associated to $X$ by the co-adjoint action of $U$:
$X^{\dagger}_{\zeta}=X^{\dagger}_{\mu(x)}=T_x\mu.X^{\#}_x=0$, so that $X\in
\mathfrak{u}_{\zeta}=Lie(U_{\zeta})$. We have thus proved that $T_x p.v \in \ker
\wred_{[x]}$ if and only if $v \in \{X^{\#}_x\, :\, X\in \mathfrak{u}_{\zeta}\}$.
Consequently, for such a $v$, one has $T_x p.v=0$,
so that $\wred$ is non-degenerate and $\mu^{-1}(\{\zeta\})/U_{\zeta}$ is a
symplectic manifold. When $\zeta=0\in\mathfrak{u}^*$, $U_{\zeta}=U$ and one
usually denotes $\mu^{-1}(\{0\})/U$ by $M\quot U$. This manifold is called the
\emph{symplectic quotient} of $M$ by $U$. Observe that in this case
$\mu^{-1}(\{0\})$ is a co-isotropic submanifold of $M$, since, if $\mu(x)=0$,
then for all $X\in\mathfrak{u}$, $T_x\mu.X^{\#}_x=X^{\dagger}_0 = 0$, so
that $(\ker T_x\mu)^{\perp_{\w}}=T_x(U.x) \subset \ker T_x\mu$. And the $2$-form
$\wred$ is then symplectic because the leaves of the null-foliation of
$\w|_N$ (that is, the foliation corresponding to the distribution $x\mapsto
\ker (\w|_N)_x = (T_x N)^{\perp_{\w}}= (\ker T_x \mu)^{\perp_{\w}}$) are precisely the $U$-orbits.\vskip 9pt

In \cite{LS}, the authors study the case where regularity
assumptions (such as assuming the action of $U$ on
$\mu^{-1}(\{0\})$ to be free, or the weaker assumption that $0$ is a regular value of $\mu$) are dropped. More precisely, Lerman and
Sjamaar showed that when
the above regularity assumptions are dropped, the reduced space
$M//U$ is a
union of symplectic manifolds which are the strata of a stratified
space. Their proof relies on the Guillemin-Marle-Sternberg normal form for the momentum
map. See subsection \ref{stratif} for further comments.

\subsection{The smooth case}

Let us now come back to the quasi-Hamiltonian setting. In \cite{AMM}, Alekseev,
Malkin and Meinrenken showed how to construct new quasi-Hamiltonian spaces
from a given quasi-Hamiltonian $U$-space $(M,\w,\mu:M\to U)$ by a reduction
procedure, assuming that $U$ is a product group $U=U_1\times U_2$ (so that
$\mu$ has two components $\mu=(\mu_1,\mu_2)$). Their result says that the
reduced space $\mu_1^{-1}(\{u\})/(U_1)_u$ is a quasi-Hamiltonian $U_2$-space. As a special case, when $U_2=\{1\}$, they obtain a \emph{symplectic manifold}.
Since this is the case we are interested in, we will state their result in
this way and give a proof that is valid in this particular situation. We
refer to \cite{AMM} for the general case. It is quite remarkable that one
can obtain symplectic manifolds from quasi-Hamiltonian spaces by a reduction
procedure. As a matter of fact, this is one of the nicest features of the
notion of quasi-Hamiltonian spaces: it enables one to obtain symplectic
structures on quotient spaces (typically, moduli spaces) using simple finite
dimensional objects as a total space. The most important example in that
respect is the moduli space of flat connections on a Riemann surface
$\Sigma$, first obtained (in the case of a compact surface) by Atiyah and
Bott in \cite{AB} by symplectic reduction of an infinite-dimensional
symplectic manifold. We refer to \cite{AMM} and to section
\ref{sympstructure} below to see how one can recover these symplectic structures
using quasi-Hamiltonian spaces. Let us now state and prove the result we are
interested in.
\begin{thm}[Symplectic reduction of quasi-Hamiltonian
spaces, the smooth case, \cite{AMM}]\label{reduction}
Let $(M,\w,\mu:M\to U)$ be a quasi-Hamiltonian $U$-space. Assume that $U$ acts freely
on $\mu^{-1}(\{1\})$. Then  $1$ is a regular value of $\mu$. Further, let
$i:\mu^{-1}(\{1\})\hookrightarrow M$ be the inclusion of the level manifold
$\mu^{-1}(\{1\})$ in $M$ and let $p:\mu^{-1}(\{1\})\to \mu^{-1}(\{1\})/U$ be
the projection on the orbit space. Then there exists a unique $2$-form $\w^{red}$ on the reduced manifold
$M^{red}:=\mu^{-1}(\{1\})/U$ such that $p^*\w^{red}=i^*\w$ on $\mu^{-1}(\{1\})$ and this $2$-form $\w^{red}$ is symplectic.
\end{thm}
We call this case the \emph{smooth case} because in this case the quotient is a smooth manifold. We see from the statement of the theorem that this case arises when the action of $U$ on $\mu^{-1}(\{1\})$ is a \emph{free action}.
\begin{proof}
Take $x\in \mu^{-1}(\{1\})$. Then, by Proposition \ref{facts}, one has $\mathrm{Im~}T_x\mu =\mathfrak{u}_x^{\perp}$. Since the action of $U$ on $\mu^{-1}(\{1\})$ is free, one has $\mathfrak{u}_x=0$ and therefore $\mathrm{Im~}T_x\mu =\mathfrak{u}$. Consequently, $1\in U$ is a regular value of $\mu$ and $\mu^{-1}(\{1\})$ is a submanifold of $M$.
The end of the proof consists in showing that $i^*\w$ is basic
with respect to the principal fibration $p$ and then verifying that the unique $2$-form
$\wred$ on $\mu^{-1}(\{1\})/U$ such that $p^*\wred=i^*\w$ is indeed
symplectic.\\ Let us first show that $i^*\w$ is basic:
\begin{eqnarray*}
& u^*(i^*\w) = i^*\w & \mathrm{for\ all\ } u\in U\\
\mathrm{and} & & \\
& \iota_{X^{\#}}i^*\w=0 & \mathrm{for\ all\ } X\in \mathfrak{u}
\end{eqnarray*}
 The first condition is obvious since $\w$ is $U$-invariant.
Consider now $X\in\mathfrak{u}$. Then:
\begin{eqnarray*}
\iota_{X^{\#}}(i^*\w) & = & i^*(\iota_{X^{\#}}\w) \\
& = & i^* \big(\frac{1}{2} \mu^*(\theta^L + \theta^R \,|\, X)\big)\\
& = & \frac{1}{2} \big( i^*\circ \mu^* (\theta^L + \theta^R\,|\, X)\big) \\
& = & \frac{1}{2} (\mu\circ i)^* (\theta^L + \theta^R \,|\, X)\\
& = & 0
\end{eqnarray*}
 since $\mu\circ i$ is constant on $\mu^{-1}(\{1\})$ and therefore
$T(\mu\circ i)=0$, hence $(\mu\circ i)^*=0$. Then there exists, by
Lemma \ref{basicforms}, a unique $2$-form $\wred$ on $\mu^{-1}(\{1\})/U$
such that $p^*\wred=i^*\w$.\\
Let us now prove that $\wred$ is a symplectic form. First:
\begin{eqnarray*}
p^*(d\wred) & = & d(p^*\wred) \\
& = & d(i^*\w) \\
& = & i^*(d\w) \\
& = & i^*(-\mu^*\chi) \\
& = & -\underset{=0}{\underbrace{(\mu \circ i)^*}} \ \chi \\
& = & 0
\end{eqnarray*}
 so that $d\wred=0$. Second, take $[x]\in \mu^{-1}(\{1\})/U$, where
$x\in \mu^{-1}(\{1\})$, and $[v]\in \ker \wred_{[x]}$, where $v\in
T_x\mu^{-1}(\{1\})=\ker T_x \mu$. Then, for all $w\in
T_x\mu^{-1}(\{1\})=\ker T_x\mu$, one has: $$(i^*\w)_x(v,w)= (p^*\wred)_x(v,w)=
\wred_{[x]}([v],[w])= 0$$ since $[v] \in \ker \wred_{[x]}$. 
Hence:
\begin{eqnarray*}
v \in \ker (i^*\w)_x & = & \{ s \in \ker T_x \mu \ | \  \forall \, w\in \ker
T_x\mu,\ \w_x(s,w)=0 \}\\
& = & \ker T_x\mu \cap (\ker T_x\mu)^{\perp_{\w}} \subset T_x M
\end{eqnarray*}
 But, by Proposition \ref{facts}, $(\ker T_x\mu)^{\perp_{\w}}= \{X^{\#}_x\, :\,
X\in\mathfrak{u}\}$, so $v=X^{\#}_x$ for some $X\in\mathfrak{u}$. Hence:
$$[v]= T_x p.v = T_x p.X^{\#}_x = 0$$ so that $\wred$ is non-degenerate.
\end{proof}

\subsection{The stratified case}\label{stratif}

What happens if we now drop the regularity assumptions of Theorem
\ref{reduction}? First one may observe that if instead of assuming the action of $U$ on $\mu^{-1}(\{1\})$ to be free one assumes that $1$ is a regular value of $\mu$, then one still has $\mathfrak{u}_x=(\mathrm{Im~} T_x\mu)^{\perp}=\{0\}$ so that the stabilizer $U_x$ of any $x\in \mu^{-1}(\{1\})$ is a discrete, hence finite (since $U$ is compact), subgroup of $U$. Consequently, $\mu^{-1}(\{1\})/U$ is a symplectic orbifold (this is the point of view adopted in \cite{AMM}). Following the techniques used in \cite{LS}
for usual Hamiltonian spaces, we will show that if we do not assume that $U$ acts freely on
$\mu^{-1}(\{1\})$ nor that $1$ is
a regular value of $\mu:M\to U$ then the orbit space $\mu^{-1}(\{1\})/U$ is a disjoint
union, over subgroups $K\subset U$, of symplectic manifolds $(N_K')^{red}$:
$$\mu^{-1}(\{1\})/U=\bigsqcup_{K\subset U} (N_K')^{red}$$ each $(N_K')^{red}$ being
obtained by applying Theorem \ref{reduction} to a quasi-Hamiltonian
space $(N_K',\w_K,\widehat{\mu_K}':N_K'\to L_K)$. Actually, the study conducted in
\cite{LS} is far more precise and ensures that the reduced space
$M^{red}:=\mu^{-1}(\{1\})/U$ is a \emph{stratified space} $M^{red}=\cup_{K\subset U} S_K$(in particular,
there is a notion of \emph{smooth function} on $M^{red}$, and the set
$\calC^{\infty}(M^{red})$ of smooth functions is an algebra over the field $\R$), with strata $(S_K)_{K\subset U}$, such that:
\begin{enumerate}
\item[-] each stratum $S_K$ is a symplectic manifold (in particular $\calC^{\infty}(S_K)$ is a Poisson algebra).
\item[-] $\calC^{\infty}(M^{red})$ is a Poisson algebra.
\item[-] the restriction maps $\calC^{\infty}(M^{red})\to \calC^{\infty}(S_K)$ are Poisson maps.
\end{enumerate}
 A stratified space satisfying these additional three conditions is
called a \emph{stratified symplectic space}. In \cite{LS}, to show that
$M^{red}$ is always a stratified symplectic space, Lerman and Sjamaar actually
obtain this space as a disjoint union of symplectic manifolds in two
differents ways. The first one enhances the stratified structure of
$M^{red}$ (the stratification being induced by the partition of $M$
according to orbit types for the action of $U$), and relies on the Guillemin-Marle-Sternberg normal form
for the momentum map. It also shows that each stratum carries a symplectic
structure. The second description of $M^{red}$ as a disjoint union of
symplectic manifolds then aims at relating this reduction to the regular
Marsden-Meyer-Weinstein procedure: the symplectic structure on each stratum is obtained by \emph{symplectic reduction} from a submanifold of $M$ endowed with a \emph{free} action of a compact Lie group. We also refer to \cite{OR} for a detailed account on the stratified symplectic structure of symplectic quotients in usual Hamiltonian geometry.\vskip 9pt
Here, we shall not be dealing with the notion of stratified space and we
will content ourselves with a description of $\mu^{-1}(\{1\})/U$ as a
disjoint union of symplectic manifolds obtained \emph{by reduction} from a
quasi-Hamiltonian space $N_K'\subset M$. We will nonetheless call the case at
hand the stratified case.

\subsubsection{\textbf{Isotropy submanifolds}}\label{isotropy_submanifolds}

We start with a quasi-Hamiltonian space $(M,\w,\mu:M\to U)$ and use the partition of $M$ given by what we may call
the \emph{isotropy type}: $$M=\bigsqcup_{K\subset U} M_K$$ where $K\subset
U$ is a closed subgroup of $U$ and $M_K$ is the set of points of $M$ whose
stabilizer is exactly $K$: $$M_K=\{x\in M \ |\ U_x=K\}.$$ Observe that if
one wants $K$ to be the stabilizer of some $x\in M$, one has to assume that
$K$ is closed, since a stabilizer always is. If $M_K$ is
non-empty, it is a submanifold of $M$ (see Proposition
\cite{GS}, p.203), called the \emph{manifold of symmetry} $K$ in
\cite{LS}. As for us, we will follow \cite{OR} and call $M_K$ the \emph{isotropy submanifold of type} $K$. The tangent space at some point $x\in M_K$ consists of all
vectors in $T_x M$ which are fixed by $K$:
$$T_x M_K = \{v\in T_x M\ |\ \mathrm{for\ all}\ k\in K, k.v=v\}$$
where $k\in K$ acts on $T_x M$ as the tangent map of the diffeomorphism $y\in
M\mapsto k.y$ which sends $x$ to itself by definition.
The action of $U$ does not preserve $M_K$ but $M_K$ is globally
stable under the action of elements $n\in \mathcal{N}(K)\subset U$, where
$\mathcal{N}(K)$ denotes the normalizer of $K$ in $U$:
$$\mathcal{N}(K):=\{u\in U \ |\ \mathrm{for\ all}\ k\in K, uku^{-1} \in
K\}.$$ It is actually the largest subgroup of $U$ leaving $M_K$ invariant,
since the stabilizer of $u.x$ for some $x\in M_K$ and some $u\in U$ is still
$U_x$ if and only if $uU_xu^{-1}=U_x$, that is, $uKu^{-1}=K$. Observe that
we have: $$Lie\big(\mathcal{N}(K)\big)\subset\{X \in\mathfrak{u}\ |\ \mathrm{for\ all}\ Y \in
\mathfrak{k}, [X,Y]\in\mathfrak{k}\}.$$ That is, the Lie algebra of the
normalizer of $K$ in $U$ is included in the normalizer of $n(\mathfrak{k})$ of the
Lie algebra $\mathfrak{k}:=Lie(K)$ in $\mathfrak{u}=Lie(U)$. The subgroup
$K$ is normal in $\mathcal{N}(K)$ and acts trivially on $M_K$ by definition
of the isotropy submanifold of type $K$, so that $M_K$ inherits an action of the
quotient group $\mathcal{N}(K)/K$. It actually follows from the definition
of $M_K$ that this induced action is free: if $n\in\mathcal{N}(K)$
stabilizes some $x$ in $M_K$, then $n\in K$ and so is the identity in
$\mathcal{N}(K)/K$. We now wish to show that $M_K$ is a quasi-Hamiltonian
space with respect to this action. We need to find a momentum map
$\mu_K:M_K\to \mathcal{N}(K)/K$ and a $2$-form $\w_K$ satisfying the axioms
of definition \ref{def-qHam}. The natural candidates are $\mu_K:=\mu|_{M_K}$
and $\w_K:=\w|_{M_K}$, but the problem is that $\mu_K$ does
not take its values in $ \mathcal{N}(K)/K$. We will now show that 
$\mu(M_K)\subset \mathcal{N}(K)$ and that we can therefore consider the
composed map $\widehat{\mu_K}:=p_K\circ\mu_K:M_K\to\mathcal{N}(K)/K$, where
$p_K$ is the projection map $p_K:\mathcal{N}(K)\to\mathcal{N}(K)/K$. Denote
then by $L_K$ the group $L_K:=\mathcal{N}(K)/K$. As $K$ is closed in $U$, so
is $\mathcal{N}(K)$, and since $U$ is compact, $\mathcal{N}(K)$ is compact.
Therefore $L_K=\mathcal{N}(K)/K$ is a compact Lie group. We will then show
that $(M_K,\w|_{M_K},\widehat{\mu_K})$ is a quasi-Hamiltonian space. Moreover, we
will show that $1\in L_K$ is a regular value of $\widehat{\mu_K}$ and that
$L_K$ acts freely on $\widehat{\mu_K}^{-1}(\{1\})$, so
that, by Theorem \ref{reduction}, the reduced space
$M_K^{red}:=\widehat{\mu_K}^{-1}(\{1\})/L_K$ is a symplectic manifold.\vskip 9pt To do
so, we start by studying $\mu(M_K)$. This whole analysis adapts the ideas of
\cite{LS} to the quasi-Hamiltonian setting. Let us denote $\w_K:=\w|_{M_K}$
and $\mu_K:=\mu|_{M_K}$. First, for all
$X\in\mathfrak{k}$, we have:
\begin{eqnarray}\label{momentumforL}
& \iota_{X^{\#}}\w_K = \frac{1}{2} \mu_K^* (\theta^L + \theta^R \, |\,X) &
\end{eqnarray}
 (where $\theta^L$ and $\theta^R$ denote as usual the Maurer-Cartan
$1$-forms of $U$, so that the above
relationship simply follows from the fact that $(M,\w,\mu:M\to U)$ is a
quasi-Hamiltonian space). Second, since $K$ acts trivially on $M_K$, we have,
for all $x\in M_K$ and all $k\in K$:
$$\mu_K(x)=\mu_K(k.x)=k\mu_K(x)k^{-1}$$ so that $\mu(x)$ belongs to the
centralizer of $K$ in $U$:
$$\mathcal{C}(K):=\{u
\in U\ |\ \mathrm{for\ all}\ k\in K, uku^{-1}=k\} $$
Since $\mathcal{C}(K)\subset\mathcal{N}(K)$, we have: $$\mu(M_K)\subset \mathcal{C}(K)\subset \calN(K).$$ We can therefore consider the map
$\widehat{\mu_K}:=p_K\circ\mu_K:M_K\to L_K=\mathcal{N}(K)/K$, where
$p_K:\mathcal{N}(K)\to\mathcal{N}(K)/K$. Furthermore, we may identify the
Lie algebra of $L_K$ to $Lie(\mathcal{N}(K))/\mathfrak{k}$. Under this
identification, the Maurer-Cartan $1$-forms $\theta^L_{L_K}$ and
$\theta^R_{L_K}$ of $L_K$ are obtained by restricting those of $U$ to
$\mathcal{N}(K)$ (which gives $Lie(\mathcal{N}(K))$-valued $1$-forms) and
composing by the projection $Lie(\mathcal{N}(K))\to
Lie(\mathcal{N}(K))/\mathfrak{k}$. It is then immediate from relation
(\ref{momentumforL}), that for all $X\in Lie(L_K)$, one has:
\begin{eqnarray}\label{momentum_axiom}
\iota_{X^\#}\w_K=\frac{1}{2}\widehat{\mu_K}^*(\theta^L_{L_K}+\theta^R_{L_K}\
|\ X)\end{eqnarray} Likewise, the Cartan $3$-form $\chi_{L_K}$ of $L_K$ is obtained by
restricting that of $U$ to $\mathcal{N}(K)$ and composing the
$Lie(\mathcal{N}(K))$-valued $3$-form thus obtained by the projection
$Lie(\mathcal{N}(K))\to Lie(\mathcal{N}(K))/\mathfrak{k}$. Then, it follows from the
fact that $d\w=-\mu^*\chi$ that we have:
\begin{eqnarray}\label{differential_axiom}
d\w_K=-\mu_K^*\chi|_{\mathcal{N}(K)}=-\widehat{\mu_K}^*\chi_{L_K}\end{eqnarray}
Thus, we have almost proved that $(M_K,\w_K,\widehat{\mu_K})$ is a
quasi-Hamiltonian $L_K$-space. In order to compute $\ker (\w_K)_x$ for all
$x\in M_K$, we observe the following two facts, the first of which is
classical in symplectic geometry and the second of which is a
quasi-Hamiltonian analogue:
\begin{lem}\label{cassymp}
Let $(V,\w)$ be a symplectic vector space and let $K$ be a compact group
acting linearly on $V$ preserving $\w$. Then the subspace $$V_K:=\{v\in V\
|\ \mathrm{for\ all}\ k\in K, k.v=v\}$$ of $K$-fixed vectors in $V$ is a
symplectic subspace of $V$.
\end{lem}
\begin{proof}
Since $K$ is compact, there exists a $K$-invariant positive definite scalar
product on $V$, that we shall denote by $(.\, |\, .)$. Since $\w$ is
non-degenerate, there exists, for any $v\in V$, a unique vector $Av\in V$
satisfying $$(v\, |\,w)=\w(Av,w)$$ for all $w\in V$, and the map $A:V\to V$
thus defined is an automorphism of $V$. Moreover, it satisfies
$A(V_K)\subset V_K$. Indeed, if $v\in V_K$, then for all $k\in K$, one has,
for  all $w\in V$:
\begin{eqnarray*}
\w(k.Av,w) & = & \w(Av,k^{-1}.w) \\
& = & (v\, |\,k^{-1}.w) \\
& = & (k.v\, |\,w) \\
& = & \w(A(k.v),w) \\
& = & \w(Av,w)
\end{eqnarray*}
 and therefore $k.Av=Av$ for all $k\in K$ (incidentally, if one
forgets the last equality, which used the fact that $k.v=v$, this also proves
that $Ak=kA$ for all $k\in K$), hence $Av\in V_K$. If now $v\in V_K$
satisfies $\w(v,w)=0$ for all $w\in V_K$, then in particular for $w=Av$, one
obtains $\w(v,Av)=0$, that is, $(v\, |\,v)=0$, hence $v=0$, since $(.\, |\, .)$
is positive definite.
\end{proof}
\begin{lem}\label{casqHam}
Let $(V,\w)$ be a vector space endowed with a possibly degenerate
antisymmetric bilinear form and let $K$ be a compact group acting linearly
on $V$ preserving $w$. Then the $2$-form $w_K:=\w|_{V_K}$ defined on the
subspace $$V_K:=\{v\in V\ |\ \mathrm{for\ all}\ k\in K, k.v=v\}$$ of
$K$-fixed vectors of $V$ has kernel: $$\ker \w_K = \ker \w \cap V_K$$ 
\end{lem}
\begin{proof}
If $\w$ is non-degenerate then this is simply Lemma \ref{cassymp}. Assume
now that $\ker \w\not= \{0\}$. Observe that $\ker \w_K=V_K^{\perp_{\w}}\cap
V_K \supset \ker \w \cap V_K$. We now consider the reduced vector space
$V^{red}:=V/\ker\w$. The $2$-form $\w$ induces a $2$-form $\w^{red}$ on
$V^{red}$, which is non-degenerate by construction. The map
$V_K\hookrightarrow V\to V/\ker\w$ induces an inclusion $V_K/(\ker\w\cap V_K)
\hookrightarrow V/\ker\w$. Further, the action of $K$ on $V$ induces an
action $k.[v]:=[k.v]$ on $V^{red}$: this action is well-defined because $K$
preserves $\w$ and therefore if $r\in \ker\w$ then $k.r\in \ker\w$. The
subspace $(V^{red})_K$ of $K$-fixed vectors for this action can be
identified with $V_K/(\ker\w\cap V_K)$. Indeed, if $[v]\in V^{red}$
satisfies, for all $k\in K$, $[k.v]=[v]$, then set:
$$w:=\int_{k\in K}(k.v)d\lambda(k)$$ where $\lambda$ is the Haar measure on
the compact Lie group $K$ (such that $\lambda(K)=1$). Then for all $k'\in K$:
\begin{eqnarray*}
k'.w & = & k'.\Big(\int_{k\in K} (k.v)d\lambda(k)\Big) \\
& = & \int_{k\in K} (k'k.v)d\lambda(k) \\
& = &  \int_{h\in K} (h.v) d\lambda(h) \\
& = & w
\end{eqnarray*}
 since the Haar measure on $K$ is invariant by translation. Thus
$w\in V_K$ and we have:
\begin{eqnarray*}
[w] & = & \Big[ \int_{k\in K} (k.v) d\lambda(k)\Big] \\
& = & \int_{k\in K} \underset{=[v]}{\underbrace{[k.v]}} d\lambda(k)\\
& = & [v] \times \int_{k\in K}d\lambda(k) \\
& = & [v].
\end{eqnarray*}
 Thus $[v] \in V_K / (\ker \w \cap V_K) \subset V^{red}$, which
proves that $(V^{red})_K \subset V_K / (\ker\w \cap V_K)$, and
therefore: $$(V^{red})_K = V_K / (\ker\w \cap V_K)$$ (the converse inclusion
being obvious). Consequently, since $V^{red}$ is a symplectic space, Lemma
\ref{cassymp} applies and we obtain: $$\ker \w^{red}|_{(V^{red})_K} =\{0\}.$$
Now $\w_K=\w|_{V_K}$ induces a $2$-form $(\w_K)^{red}$ on $V_K/(\ker\w \cap
V_K)=(V^{red})_K$, whose kernel is, by definition: $$\ker (\w_K)^{red} = \ker \w_K / (\ker
\w \cap V_K).$$ But, again by definition, $(\w_K)^{red} =
\w^{red}|_{(V^{red})_K}$, so that $\ker (\w_K)^{red} = \{0\}$, hence
$\ker\w_K = \ker\w \cap V_K$, which proves the lemma. 
\end{proof}
We then obtain a new class of examples of quasi-Hamiltonian spaces:
\begin{thm}\label{strata_one}
For each closed subgroup $K\subset U$, the compact Lie group
$L_K:=\mathcal{N}(K)/K$ acts freely on the isotropy submanifold $$M_K=\{x\in
M\ |\ U_x =K\}.$$ In addition to that, $\mu(M_K)\subset \mathcal{N}(K)$ and
$(M_K,\w_K:=\w|_{M_K},\widehat{\mu_K}:=p_K\circ\mu|_{M_K})$, where $p_K$ is
the projection map $p_K:\mathcal{N}(K)\to\mathcal{N}(K)/K=L_K$, is a
quasi-Hamiltonian space.
\end{thm}
\begin{proof}
Observe first that $\widehat{\mu_K}$ is $L_K$ equivariant because $\mu$ is $U$-equivariant and $p_K:\calN(K)\to\calN(K)/K$ is a group morphism. Second, recall that we have obtained the relations (\ref{momentum_axiom}) and (\ref{differential_axiom}), so that, to prove that $(M_K,\w_K,\widehat{\mu_K}:M_K\to L_K)$ is a quasi-Hamiltonian $L_K$-space,
the only thing left to do is compute $\ker (\w_K)_x\subset T_xM_K$. Since $T_xM_K=\{x\in T_xM\ |\ \forall k\in K, k.v=v\}$, Lemma \ref{casqHam} applies and one has: 
$$\ker (\w_K)_x = \ker \w_x \cap T_x M_K = \{X^{\#}_x\, :\, X\in \mathfrak{u}\ |\ Ad\, \mu(x).X=-X\} \cap T_x M_K.$$ But a vector of $T_x M$ of the form $X^{\#}_x$ lies in $T_x M_K\subset T_x M$ if and only of $X\in Lie(\mathcal{N}(K))\subset \mathfrak{u}$. Further, we have seen that for all $x\in M_K$, $\mu(X)=\mu_K(x)\in \mathcal {N}(K)$. Therefore: $$\ker (\w_K)_x = \{X^{\#}_x\, :\, X\in Lie(\mathcal{N}(K))\ |\ Ad\, \mu_K(x).X=-X\}.$$ Since $K$ acts trivially on $M_K$ and on $\mathcal{N}(K)/K$, this last statement is equivalent to: $$\ker (\w_K)_x =
\{X^{\#}_x\, :\,X\in Lie(\mathcal{N}(K))/\mathfrak{k} \ |\
Ad\,\widehat{\mu_K}(x).X = -X\}$$ 
which completes the proof.
\end{proof}
And we then observe that:
\begin{cor}\label{strata_two}
$1\in L_K$ is a regular value of $\widehat{\mu_K}$ and the
reduced space $M_K^{red}:=\widehat{\mu_K}^{-1}(\{1\})/L_K$ is a symplectic
manifold.
\end{cor}
\begin{proof}
Since the action of $L_K$ on $M_K$ is free, the fact that
$M_K^{red}:=\widehat{\mu_K}^{-1}(\{1\})/L_K$ is a symplectic
manifold follows from Theorem \ref{reduction}.
\end{proof}

\subsubsection{\textbf{Structure of quasi-Hamiltonian quotients}}

We will now use the above analysis to show that, without any regularity assumptions on the action of $U$ on $M$ or on the momentum map $\mu:M\to U$, the orbit space $M^{red}:=\mu^{-1}(\{1\})/U$ is a disjoint union of symplectic manifolds.
First, in analogy with \cite{LS}, we observe:
\begin{lem}\label{quotient_description_one}
Denote by $(K_j)_{j\in J}$ a system of representatives of conjugacy classes of closed subgroups of $U$ (every closed subgroup $K\subset U$ is conjugate to exactly one of the pairwise non-conjugate $K_j$). Denote by $M_{K_j}$ the isotropy submanifold of type $K_j$ in the quasi-Hamiltonian space $(M,\w,\mu:M\to U)$: $$M_{K_j}=\{x\in M\ |\ U_x=K_j\}.$$
Then, the orbit space $\mu^{-1}(\{1_U\})/U$ is the disjoint union:
$$\mu^{-1}(\{1_U\})/U=\bigsqcup_{j\in J} (\mu^{-1}(\{1_U\})\cap U.M_{K_j})/U.$$
\end{lem}
\begin{proof} Take a $U$-orbit $U.x$ in $\mu^{-1}(\{1_U\})$. The stabilizer $U_x$ of $x$ is conjugate to one of the $(K_j)$, that is: $U_x=uK_ju^{-1}$ for some $u\in U$. Therefore, the stabilizer of $y:=u^{-1}.x\in \mu^{-1}(\{1_U\})$ is exactly $K_j$, and we then have $U.y=U.x$ with $y\in M_{K_j}$. Therefore, we have shown: $$\mu^{-1}(\{1_U\})/U=\bigcup_{j\in J} (\mu^{-1}(\{1_U\})\cap U.M_{K_j})/U.$$
The above union is disjoint because if $U.x$ is a $U$-orbit in $\mu^{-1}(\{1_U\})\cap U.M_{K_j}$, the stabilizer of $x$ is conjugate to $K_j$ and therefore not conjugate to any $K_{j'}$ for $j'\not= j$.
\end{proof}
We will now study each one of the sets $(\mu^{-1}(\{1_U\})\cap U.M_{K_j})/U$ separately. We will show, in analogy with the result of Lerman and Sjamaar in \cite{LS}, that each one of these sets is a smooth manifold that carries a symplectic structure, and that this symplectic structure may be obtained \emph{by reduction} from a quasi-Hamiltonian space endowed with a \emph{free} action of a compact Lie group (that is, by applying Theorem \ref{reduction}). In \cite{OR}, this procedure is called \emph{Sjamaar's principle}. The way this principle is developped in \cite{OR} is way more general than what we do here: they consider the quotients $\mu^{-1}(\{\xi\})/U_{\xi}$ for an arbitrary $\xi\in\mathfrak{u}^*$, which also makes the situation slightly more complicated (notably to find an \emph{equivariant} momentum map for the isotropy submanifolds $M_K$). Here, we we begin by observing the following fact:
\begin{lem}\label{quotient_description_two}
Let $K\subset U$ be a closed subgroup of $U$ and denote by $M_K$ the isotropy submanifold of type $K$ in the quasi-Hamiltonian space $(M,\w,\mu:M\to U)$: $$M_K=\{x\in M\ |\ U_x=K\}.$$ Denote by $\calN(K)$ the normalizer of $K$ in $U$ and by $L_K$ the quotient group $L_K=\calN(K)/K$. Then, the map:
\begin{eqnarray*}
f_K: (\mu^{-1}(\{1_U\})\cap M_K)/L_K & \longto & (\mu^{-1}(\{1_U\})\cap U.M_K)/U\\
L_K.x & \longmapsto & U.x
\end{eqnarray*}
sending the $L_K$-orbit of a point $x\in (\mu^{-1}(\{1_U\})\cap M_K)$ to its $U$-orbit in $(\mu^{-1}(\{1_U\})\cap U.M_K)$ is well-defined and is a bijection: $$(\mu^{-1}(\{1_U\})\cap M_K)/L_K \overset{\simeq}{\longto}(\mu^{-1}(\{1_U\})\cap U.M_K)/U$$ Consequently, we deduce from Lemma \ref{quotient_description_one} that: $$\mu^{-1}(\{1_U\})/U=\bigsqcup_{j\in J} (\mu^{-1}(\{1_U\})\cap M_{K_j})/L_{K_j}.$$
\end{lem}
\begin{proof} The map $f_K$ is well-defined because if $x,y\in \mu^{-1}(\{1_U\})\cap M_K$ lie in a same $L_K$-orbit then they lie in a same $U$-orbit in $(\mu^{-1}(\{1_U\})\cap U.M_K)$.\\ The map $f_K$ is onto because a $U$-orbit in $(\mu^{-1}(\{1_U\})\cap U.M_K)$ is of the form $U.x$  for some $x\in (\mu^{-1}(\{1_U\})\cap M_K)$, and $f_K$ then sends the $L_K$-orbit of such an $x$ in $(\mu^{-1}(\{1_U\})\cap M_K)$ to the $U$-orbit $U.x$ in $(\mu^{-1}(\{1_U\})\cap U.M_K)$.\\ The map $f_K$ is one-to-one because if $x,y\in (\mu^{-1}(\{1_U\})\cap M_K)$ lie in a same $U$-orbit in $(\mu^{-1}(\{1_U\})\cap U.M_K)$, say $y=u.x$ for some $u\in U$, then the stabilizer of $y$ in $(\mu^{-1}(\{1_U\})\cap U.M_K)$ is $U_y=uU_xu^{-1}$. But since $x,y\in M_K$ we have $U_x=U_y=K$, hence $u\in\calN(K)$ and $L_K.y=L_K.x$. The rest of the Proposition follows from Lemma \ref{quotient_description_one}.
\end{proof}
We will now prove that each of the sets $(\mu^{-1}(\{1_U\})\cap U.M_K)/U=(\mu^{-1}(\{1_U\})\cap M_K)/L_K$ is a \emph{smooth}, \emph{symplectic manifold}. To do so, we will show that each of these sets is the quasi-Hamiltonian quotient $N'_K//L_K$ associated to a quasi-Hamiltonian space of the form $(N'_K,\w_K,\widehat{\mu_K}':N_K'\to L_K)$ (see Theorem \ref{strata_one} and Corollary \ref{strata_two}). More precisely, we have to show that $$(\mu^{-1}(\{1_U\})\cap M_K)/L_K = (\widehat{\mu_K}')^{-1}(\{1_{L_K}\})/L_K$$ where $\widehat{\mu_K}'$ is the momentum map of a \emph{free} action of $L_K$ on a quasi-Hamiltonian space $(N_K',\w_K,\widehat{\mu_K}':N_K'\to L_K)$. This last step is not entirely immediate. In fact, experience from the usual Hamiltonian case dealt with by Lerman and Sjamaar in \cite{LS} shows that in that setting too, one has to replace $(M_K,\w_K,\widehat{\mu_K}:M_K\to Lie(L_K)^*)$ by another Hamiltonian $L_K$-space $(M_K',\w_K,\widehat{\mu_K}':M_K'\to Lie(L_K)^*)$, that space $M_K'$ being the union of connected components of $M_K$ which have a non-empty intersection with $\mu^{-1}(\{0\})$. The point is that this space $M_K'$ is in a way big enough to study the quotient $(\mu^{-1}(\{0\})\cap M_K)/L_K$ because by definition of $M_K'$ one has $(\mu^{-1}(\{0\})\cap M_K)/L_K=(\mu^{-1}(\{0\})\cap M_K')/L_K$. And then one can prove that $(\mu^{-1}(\{0\})\cap M_K')/L_K=\widehat{\mu_K}'^{-1}(\{0\})/L_K=(M_K')^{red}$ (whereas it is not true that $(\mu^{-1}(\{0\})\cap M_K)/L_K=\widehat{\mu_K}^{-1}(\{0\})/L_K$), thereby proving that $(\mu^{-1}(\{0\})\cap M_K)/L_K=(M_K')^{red}$ is a symplectic manifold. Trying an exactly analogous approach in the quasi-Hamiltonian setting does not work: the union of connected components of $M_K$ containing points of $\mu^{-1}(\{1_U\})$ is still too big, and one has to introduce another quasi-Hamiltonian $L_K$-space, which we will denote by $N_K$ (see Lemma \ref{same_quotient}). This is what we do next (see also remark \ref{small_enough}). We begin with the following lemma:
\begin{lem}\label{smaller}
Let $\calB\subset \mathfrak{u}$ be an $Ad$-stable open ball centered at $0\in\mathfrak{u}$ such that the exponential map $\exp|_{\calB}:\calB\to \exp(\calB)$ is a diffeomorphism onto an open subset of $U$ containing $1_U$. Denote by $N\subset M$ the $U$-stable open subset of $M$ defined by $$N:=\mu^{-1}(\exp(\calB)).$$ Then $(N,\w|_N,\mu|_N:N\to U)$ is a quasi-Hamiltonian $U$-space, and one has: $$(\mu|_N)^{-1}(\{1_U\})/U=\mu^{-1}(\{1_U\})/U.$$
\end{lem}
\begin{proof}
Any $U$-stable \emph{open} subset of a quasi-Hamiltonian space is a quasi-Hamiltonian space when endowed with the restriction of the $2$-form and the restriction of the momentum map. In the above case, one has $(\mu|_N)^{-1}(\{1_U\})=\mu^{-1}(\{1_U\})$ by construction of $N=\mu^{-1}(\exp(\calB))$.
\end{proof}
We can then compare the isotropy submanifolds of $M$ and of $N$:
\begin{lem}\label{same_quotient}
Let $(N,\w|_N,\mu|_N:N\to U)$ be the quasi-Hamiltonian $U$-space introduced in Lemma \ref{smaller}. Let $K\subset U$ be a closed subgroup of $U$ and denote by $$M_K=\{x\in M\ |\ U_x=K\}\mathrm{\  and\ } N_K=\{x\in N\ |\ U_x=K\}$$ the isotropy submanifolds of type $K$ of $M$ and $N$ respectively. Then one has: $$\mu^{-1}(\{1_U\})\cap M_K = \mu^{-1}(\{1_U\})\cap N_K.$$
\end{lem}
\begin{proof}
The equality $\mu^{-1}(\{1_U\})\cap M_K = \mu^{-1}(\{1_U\})\cap N_K$ follows from the fact that $\mu^{-1}(\{1_U\})\subset N$ by construction of $N=\mu^{-1}(\exp(\calB))$.
\end{proof}
We will now show that the orbit space $(\mu^{-1}(\{1_U\})\cap N_K)/L_K$ has a symplectic structure. To do this, we apply Theorem \ref{linearization} to the quasi-Hamiltonian space $M_0=N=\mu^{-1}(\exp(\calB))$ constructed in Lemma \ref{smaller} to obtain the following result:
\begin{lem}\label{final_step}
Let $(N=\mu^{-1}(\exp(\calB)),\w|_N,\mu|_N:N\to U)$ be the quasi-Hamiltonian $U$-space introduced in Lemma \ref{smaller}. Let $K\subset U$ be a closed subgroup of $U$ and let $\calN(K)$ be its normalizer in $U$. Denote by $L_K$ the quotient group $L_K:=\calN(K)/K$ and by $p_K$ the projection $p_K:\calN(K)\to L_K=\calN(K)/K$. Let $$N_K=\{x\in N\ |\ U_x=K\}$$ be the istotropy submanifold of type $K$ in $N$. Recall from Theorem \ref{strata_one} that $\mu(N_K)\subset \calN(K)$ and that $(N_K,\w|_{N_K},\widehat{\mu_K}=p_K\circ\mu|_{N_K}:N_K\to L_K)$ is a quasi-hamitonian $L_K$-space. Denote by $N_K'$ the union of connected components of $N_K$ which have a non-empty intersection with $\mu^{-1}(\{1_U\})$, and by $\widehat{\mu_K}'$ the restriction of $\widehat{\mu_K}$ to $N_K'$. Then: $N_K'$ is $L_K$-stable and $(N_K',\w|_{N_K'},\widehat{\mu_K}':N_K'\to L_K)$ is a quasi-Hamiltonian $L_K$-space. Furthermore, one has:  $$\mu^{-1}(\{1_U\})\cap N_K = \mu^{-1}(\{1_U\})\cap N_K' = (\widehat{\mu_K}')^{-1}(\{1_{L_K}\})$$ and consequently: $$(\mu^{-1}(\{1_U\})\cap N_K)/L_K = (\mu^{-1}(\{1_U\})\cap N_K')/L_K = (\widehat{\mu_K}')^{-1}(\{1_{L_K}\})/L_K=(N_K')^{red}$$
\end{lem}
\begin{proof}
We first show that $N_K'$ is $L_K$-stable and is a quasi-Hamiltonian $L_K$-space. If $x\in N_K'$ and $n\in\calN(K)$ then there exists, by definition of $N_K'$, a point $x_0\in\mu^{-1}(\{1_U\})\cap N_K$ which is connected to $x$ by a path $(x_t)$ in $N_K$. Then $(n.x_t)$ is a path from $(n.x_0)$ to $(n.x)$ in $N_K$. Since $\mu(n.x_0)=n\mu(x_0)n^{-1}=1_U$ and $(n.x)$ lies in the same connected component of $N_K$ as $(n.x_0)$, we have $(n.x)\in N_K'$. The fact that $(N_K',\w|_{N_K'},\widehat{\mu_K}':N_K'\to L_K)$ is a quasi-Hamiltonian space then follows from the fact that $N_K'$ is an $L_K$-stable open subset of the quasi-Hamiltonian space $(N_K,\w|_{N_K},\widehat{\mu_K}:N_K\to L_K)$.\\ Let us now prove that $\mu^{-1}(\{1_U\})\cap N_K = \mu^{-1}(\{1_U\})\cap N_K' = (\widehat{\mu_K}')^{-1}(\{1_{L_K}\})$. By definition of $N_K'$, one has $\mu^{-1}(\{1_U\})\cap N_K=\mu^{-1}(\{1_U\})\cap N_K'$. Furthermore, it is obvious that $\mu^{-1}(\{1_U\})\cap N_K' \subset (\widehat{\mu_K}')^{-1}(\{1_{L_K}\})$ since $\widehat{\mu_K}'=p_K\circ\widehat{\mu_K}|_{N_K'}$ and $p_K:\calN(K)\to \calN(K)/K$ is a group morphism. Let us now prove the converse inclusion. We begin by observing  that since the exponential map is invertible on $\calB\subset\mathfrak{u}$ and $N=\exp(\calB)$, Theorem \ref{linearization} applies: the map $\widetilde{\mu}:=\exp^{-1}\circ\mu|_N:N\to\mathfrak{u}$ is a momentum map in the usual sense for the action of $U$ on $N$ and $\mu^{-1}(\{1_U\})=\widetilde{\mu}^{-1}(\{0\})$. In particular, one has, for all $x\in N_K'$, $\mathrm{Im~} T_x\widetilde{\mu}=\mathfrak{u}_x^{\perp}=\mathfrak{k}^{\perp}$ and, since $0\in\widetilde{\mu}(N_K')$ by definition of $N_K'$, this implies $\widetilde{\mu}(N_K')\subset\mathfrak{k}^{\perp}$. Take now $x\in(\widehat{\mu_K}')^{-1}(\{1_{L_K}\})\subset N_K'$. This means that $\mu(x)\in (K\cap \mu(N_K'))\subset \exp(\calB)$, hence $\widetilde{\mu}(x)=\exp^{-1}\circ\mu(x)\in \mathfrak{k}\cap\widetilde{\mu}(N_K')\subset \mathfrak{k}\cap\mathfrak{k}^{\perp}=\{0\}$. Consequently, $\widetilde{\mu}(x)=0$ and therefore $\mu(x)=1_U$. Hence $(\widehat{\mu_K}')^{-1}(\{1_{L_K}\})\subset \mu^{-1}(\{1_U\})\cap N_K'$, which completes the proof. 
\end{proof}
\begin{rk}\label{small_enough}
Lemma \ref{final_step} is crucial in our proof of forthcoming Theorem \ref{reduction_strat}. Although our argument is similar to the one in \cite{LS}, where the usual Hamiltonian case is treated, extra difficulties arise to show that $\mu^{-1}(\{1_U\})\cap N_K'=(\widehat{\mu_K}')^{-1}(\{1_{L_K}\})$. In particular, we were unable to obtain such a statement involving $M_K$ or $M_K'$ instead of $N_K$ and $N_K'$. In the end this is not a problem because we proved that $\mu^{-1}(\{1_U\})\cap M_K=\mu^{-1}(\{1_U\})\cap N_K=\mu^{-1}(\{1_U\})\cap N_K'$ (see Lemma \ref{same_quotient}). The point of introducing $N_K$ (and then later $N_K'$) is to be able to linearize the quasi-Hamiltonian space that we are dealing with without changing the associated quotient. This idea was suggested to us by the reading of \cite{HJS}, where a description of quasi-Hamiltonian quotients as disjoint unions of symplectic manifolds is also obtained. The main difference between Theorem \ref{reduction_strat} and Theorem 2.9 in \cite{HJS} is that in our case the symplectic structure on each component of the union is obtained by reduction from a quasi-Hamiltonian space $(N_K',\w_K,\widehat{\mu_K}':N_K'\to L_K)$ endowed with a \emph{free} action of the compact Lie group $L_K$. The linearization theorem enables us to reduce the case at hand to the usual Hamiltonian case and mimic the argument in \cite{LS} (Theorem 3.5). It would be interesting to know if this detour can be avoided.
\end{rk}
\begin{thm}[Symplectic reduction of quasi-Hamiltonian spaces, the
stratified case]\label{reduction_strat}
Let $(M,\w,\mu:M\to U)$ be a quasi-Hamiltonian $U$-space. For any closed subgroup $K\subset U$, denote by $M_K$ the isotropy manifold of type $K$ in $M$: $$M_K=\{x\in M\ |\ U_x=K\}.$$ Denote by $\calN(K)$ the normalizer of $K$ in $U$ and by $L_K$ the quotient group $L_K:=\calN(K)/K$. Then the orbit space $$(\mu^{-1}(\{1_U\})\cap M_K)/L_K$$ is a smooth symplectic manifold.\\ Denote by $(K_j)_{j\in J}$ a system of representatives of closed subgroups of $U$. Then the orbit
space $M^{red}:=\mu^{-1}(\{1_U\})/U$ is the disjoint union of the following symplectic manifolds:
$$\mu^{-1}(\{1_U\})/U = \bigsqcup_{j\in J} (\mu^{-1}(\{1_U\})\cap M_{K_j})/L_{K_j}.$$
\end{thm}
\begin{proof}
By Lemmas \ref{same_quotient} and \ref{final_step}, we have: $$(\mu^{-1}(\{1_U\})\cap M_K)/L_{K}=(\mu^{-1}(\{1_U\})\cap N_{K})/L_{K}=(\mu^{-1}(\{1_U\})\cap N_{K}')/L_{K}=(N_K')^{red}$$ where the compact group $L_{K}$ acts freely on the quasi-Hamiltonian space $(N_K',\w_{K},\widehat{\mu_{K}}':N_{K}'\to L_{K})$, so that Theorem \ref{reduction} shows that $(\mu^{-1}(\{1_U\})\cap M_K)/L_K=(N_K')^{red}$ is a symplectic manifold. By Lemmas \ref{quotient_description_one} and \ref{quotient_description_two}, we then have:
$$\mu^{-1}(\{1_U\})/U = \bigsqcup_{j\in J} (\mu^{-1}(\{1_U\})\cap U.M_{K_j})/U = \bigsqcup_{j\in J} (\mu^{-1}(\{1_U\})\cap M_{K_j})/L_{K_j}.$$
\end{proof}
Observe that to prove that the set $(\mu^{-1}(\{1_U\})\cap M_K)/L_{K}$ is a smooth symplectic manifold, we found a quasi-Hamiltonian $L_K$-space $(N_K',\w_K,\widehat{\mu_K}':N_K'\to L_K)$  on which $L_K$ acts freely such that $(N_K')^{red}=(\mu^{-1}(\{1_U\})\cap M_K)/L_{K}$ and then applied quasi-Hamiltonian reduction in the smooth case (Theorem \ref{reduction}) to $N_K'$. One key step in this proof is to show that $(\widehat{\mu_K}')^{-1}(\{1_{L_K}\})/L_K=(\mu^{-1}(\{1_U\})\cap N_{K}')/L_{K}$ and it was to obtain this equality that we used the linearization Theorem \ref{linearization}. We then showed that for any quasi-Hamiltonian space $(M,\w,\mu:M\to U)$ the reduced space $M^{red}:= \mu^{-1}(\{1\})/U$ is a disjoint union of symplectic manifolds. We denote this reduced space
by $M//U$, as in the usual Hamiltonian case:
\begin{defi}[Quasi-Hamiltonian quotient]\label{q-ham_quot}
The reduced space $$M//U:=\mu^{-1}(\{1_U\})/U=\bigsqcup_{j\in J} (\mu^{-1}(\{1_U\})\cap M_{K_j})/L_{K_j}$$ associated, by means of
Theorems \ref{reduction} and
\ref{reduction_strat}, to a given quasi-Hamiltonian space $(M,\w,\mu:M\to U)$ is
called the \emph{quasi-Hamiltonian quotient} associated to $M$.
\end{defi}
\begin{rk}\label{regvalue}
Observe that when the action of $U$ on $M$ is free, then the only subgroup $K\subset U$ such that the isotropy submanifold $M_K$ is non-empty is $K=\{1\}$, so that the results of
Theorems \ref{reduction} and \ref{reduction_strat} do coincide in this
case.
\end{rk}
As we shall see in section \ref{sympstructure}, representation
spaces of surface groups naturally arise as quasi-Hamiltonian quotients. Since in this case it
is known that representation spaces are stratified symplectic spaces in the
sense of \cite{LS} (see for instance \cite{Hueb1}), it should be possible to
obtain this stratified symplectic structure in the quasi-Hamiltonian
framework. Following \cite{LS}, the first step to do so should be a normal
form for momentum maps on quasi-Hamiltonian spaces.

\section{Application to representation spaces of surface groups}\label{sympstructure}

In this section, we wish to briefly explain, following \cite{AMM}, how the
notion of quasi-Hamiltonian space provides a proof of the fact that,
for any Lie group $(U,(.\, |\, .))$ endowed with an $Ad$-invariant
non-degenerate product and any collection $\calC=(\calC_j)_{1\leq j\leq l}$
of $l$ conjugacy classes of $U$, there exists a symplectic structure on the
representation spaces $$\Hom_{\calC}(\pi_{g,l},U)\big/ U$$ (see
\ref{repspace} below for a precise definition of these spaces). This will serve as an example to illustrate Theorem \ref{reduction_strat}. Here, $\pi_{g,l}=\pi_1(\Sigma_{g,l})$ denotes the fundamental group of the
surface $\Sigma_{g,l}:=\Sigma_g\bs\{s_1,\,\ldots\, , s_l\}$, $\Sigma_g$
being a compact Riemann surface of genus $g\geq 0$, $l$ being an integer $l\geq 1$
and $s_1, \,{\ldots} \, ,s_l$ being $l$ pairwise distinct points of
$\Sigma_g$. When $l=0$, we set $\calC:=\emptyset$ and
$\Sigma_{g,0}:=\Sigma_g$. Everything we will say is valid for any $g\geq 0$
and any $l\geq 0$ but we will not always distinguish between the cases $l=0$
and $l\geq 1$, to lighten the presentation.\vskip 9 pt
Recall that the fundamental group of the surface
$\Sigma_{g,l}=\Sigma_g\bs\{s_1,\,{\ldots} \, ,s_l\}$ has the following
finite presentation: $$\pi_{g,l}=<\alpha_1,\,{\ldots}\, ,\alpha_g,\beta_1, \,{\ldots} \,
,\beta_g, \gamma_1,\,{\ldots} \, ,\gamma_l\ |\
\prod_{i=1}^g[\alpha_i,\beta_i]\prod_{j=1}^l\gamma_j =1>$$ each $\gamma_j$ being the
homotopy class of a loop around the puncture $s_j$. In particular, if $l\geq
1$, it is a free group on $(2g+l-1)$ generators. As a consequence of this
presentation, we see that, having chosen a set of generators of $\pi_{g,l}$, giving a \emph{representation of} $\pi_{g,l}$
\emph{in the group} $U$ (that is, a group morphism from $\pi_{g,l}$ to $U$)
amounts to giving $(2g+l)$ elements $(a_i,b_i,u_j)_{1\leq i\leq g, 1\leq
j\leq l}$ of $U$ satisfying: $$\prod_{i=1}^g[a_i,b_i]\prod_{j=1}^l u_j =1.$$ Two
representations $(a_i,b_i,u_j)_{i,j}$ and $(a'_i,b'_i,u'_j)_{i,j}$ are then
called \emph{equivalent} if there exists an element $u\in U$ such that
$a'_i=ua_iu^{-1}$, $b'_i=ub_iu^{-1}$, $u'_j=uu_ju^{-1}$ for all $i,j$. The
original approach to describing symplectic structures on spaces of
representations shows that, in order to
obtain \emph{symplectic} structures, one has to prescribe the conjugacy
class of each $u_j$, $1\leq j\leq l$. Otherwise, one may obtain Poisson
structures, but we shall not enter these considerations and refer to
\cite{Hueb3} and \cite{AKSM} instead. We are then led to studying the
space $\Hom_{\calC}(\pi_{g,l},U)$ of representations of $\pi_{g,l}$ in $U$
with prescribed conjugacy classes for the $(u_j)_{1\leq j\leq l}$:
\begin{defi}\label{repspace}
We define the space $\Hom_{\calC}(\pi_{g,l},U)$ to be the following set of
group morphisms: $$\Hom_{\calC}(\pi_{g,l},U)=\{\rho:\pi_{g,l}\to U\ |\
\rho(\gamma_j)\in\calC_j \ \mathrm{for\ all}\ j\in\{1,\,{\ldots}\, ,l\}\}.$$
\end{defi}
Observe that this space may very well be empty, depending on the
choice of the conjugacy classes $(\calC_j)_{1\leq j\leq l}$. As a matter of
fact,  when $g=0$, conditions on the $(\calC_j)$ for this set to be non-empty are quite
difficult to obtain (see for instance \cite{AW} for the case
$U=SU(n)$). However, when $g\geq 1$ and $U$ is semi-simple, the above set is always non-empty, as shown in \cite{Ho}. As earlier, giving such a morphism $\rho \in \Hom_{\calC}(\pi_{g,l},U)$ amounts to giving appropriate elements of $U$:
$$\Hom_{\calC}(\pi_{g,l},U)\simeq \{(a_1,\,{\ldots} \, ,a_g,b_1,\, {\ldots} \,
,b_g,u_1,\, {\ldots} \, u_l) \in \underset{2g\
\mathrm{times}}{\underbrace{U\times\cdots\times U}} \times \pconj \ |\
\prod_{i=1}^g[a_i,b_i]\prod_{j=1}^l u_j =1\}.$$ In particular, two
representations $(a_i,b_i,u_j)_{i,j}$ and $(a'_i,b'_i,u'_j)_{i,j}$ are
equivalent if and only if they are in a same orbit of the diagonal action of
$U$ on $U\times\cdots\times U\times\pconj$. The \emph{representation space}
$\Rep_{\calC}(\pi_{g,l},U)$ is then defined to be the quotient space for
this action: $$\Rep_{\calC}(\pi_{g,l},U):=\Hom_{\calC}(\pi_{g,l},U)\big/ U.$$ Following
for instance \cite{Hueb1}, the idea to
obtain a symplectic structure on the representation space, or moduli space,
$\Rep_{\calC}(\pi_{g,l},U)$ is then to see this quotient as a
\emph{symplectic quotient}, meaning that one wishes to identify
$\Hom_{\calC}(\pi_{g,l},U)$ with the fibre of a momentum map defined on an
\emph{extended moduli space} (the expression comes from
\cite{Jeffrey1,Hueb1}). The notion of quasi-Hamiltonian space then arises
naturally from the choice of $$\underset{2g\
\mathrm{times}}{\underbrace{U\times\cdots\times U}}\times\pconj$$ as an
extended moduli space, and of the map $$\mu_{g,l}(a_1,\, {\ldots} \,
,a_g,b_1,\,{\ldots} \, ,b_g,u_1,\,{\ldots} \, , u_l) = [a_1,b_1]{\ldots}
[a_g,b_g]u_1{\ldots} u_l$$ as $U$-valued momentum map, so that:
$$\Rep_{\calC}(\pi_{g,l},U) = \mu_{g,l}^{-1}(\{1\})/U.$$ Actually, because of
the occurrence of the commutators $[a_i,b_i]$, it is more appropriate to
re-arrange the arguments of the map $\mu_{g,l}$ in the following way:
$$\mu_{g,l}(a_1,b_1,\, {\ldots} \,
,a_g,b_g,\, ,u_1,\,{\ldots} \, , u_l) = [a_1,b_1]{\ldots}
[a_g,b_g]u_1{\ldots} u_l=1$$ and to write the extended moduli space:
$$\underset{g\
\mathrm{times}}{\underbrace{(U\times U)\cdots\times (U\times
U)}}\times\pconj.$$ In the case where $g=0$, one simply has:
$$\begin{array}{rrcl}
\mu_{0,l}: & \pconj & \longto & U \\
& (u_1,\, {\ldots} \, ,u_l) & \longmapsto & u_1{\ldots} u_l
\end{array}$$
 When $g=1$ and $l=0$, one has:
$$\begin{array}{rrcl}
\mu_{1,0}: & U\times U & \longto & U \\
& (a, b) & \longmapsto & aba^{-1}b^{-1}
\end{array}$$
 These two particular cases correspond to the examples we recalled
in Propositions \ref{pconj} and \ref{intfuseddouble}, and motivate the notion of
quasi-Hamiltonian space. Thus, in general, the extended moduli space is the
following quasi-Hamiltonian space: $$\mathcal{M}_{g,l}:=\underset{g\
\mathrm{times}}{\underbrace{\mathfrak{D}(U)\times\cdots\times\mathfrak{D}(U)}}\times\pconj.$$
(where $\mathfrak{D}(U)$ is the internally fused double of $U$ of
Proposition \ref{intfuseddouble}) equipped with the diagonal $U$-action and the momentum map
$$\begin{array}{rrcl}
\mu_{g,l}: & \mathfrak{D}(U)\times\cdots\times\mathfrak{D}(U)\times\pconj &
\longto & U \\
& (a_1,b_1,\, {\ldots} \, ,a_g,b_g,u_1, \, {\ldots} \, ,u_l) & \longmapsto &
[a_1,b_1]{\ldots}[a_g,b_g]u_1{\ldots} u_l 
\end{array}$$
 The representation space $\Rep_{\calC}(\pi_{g,l},U)$ is then the
associated quasi-Hamiltonian quotient (see definition \ref{q-ham_quot}):
$$\Rep_{\calC}(\pi_{g,l},U)=\mathcal{M}_{g,l}\quot
U=(\underset{g\
\mathrm{times}}{\underbrace{\mathfrak{D}(U)\times\cdots\times\mathfrak{D}(U)}}
\times\pconj)\quot U.$$ In particular, in the case of an $l$-punctured sphere
($g=0$),  we
have: $$\Hom_{\calC}\big(\piS,U\big)\big/ U = (\pconj)\quot U.$$ We also
spell out the case of torus: $$\Hom\big(\pi_1(\mathbb{T}^2),U\big)\big/ U =
\mathfrak{D}(U)\quot U$$ (there are no conjugacy classes necessary here, as
the surface $\mathbb{T}^2$ is closed) and of the punctured torus:
$$\Hom_{\calC}\big(\pi_1(\mathbb{T}^2\bs\{s\}\big),U)\big/ U = (\mathfrak{D}(U)\times
\calC)\quot U.$$ We then know from Theorems
\ref{reduction} and \ref{reduction_strat} that these representation spaces
$\Rep_{\calC}(\pi_{g,l},U)=\mathcal{M}_{g,l}\quot U$
carry a symplectic structure, obtained by reduction from the
quasi-Hamiltonian space $\mathcal{M}_{g,l}$. More precisely, the representation spaces $\Rep_{\calC}(\pi_{g,l},U)$ are \emph{disjoint unions of symplectic manifolds}. Observe that one essential
ingredient to obtain this symplectic structure was the fact that $\pi_{g,l}$
admits a finite presentation with a single  relation, which was used as a
momentum relation.

\end{document}